\documentclass[11pt]{amsart}

\usepackage[utf8]{inputenc}
\usepackage{mathtools}
\usepackage{amssymb, amsfonts, amsthm}
\usepackage{siunitx}
\usepackage{booktabs}
\usepackage{multirow}
\usepackage[linesnumbered,ruled,vlined]{algorithm2e}

\usepackage[
  hmarginratio={1:1},     
  textwidth=450pt,        
  heightrounded,          
]{geometry}

\usepackage{graphicx}
\usepackage{xcolor,soul}
\usepackage{url}
\usepackage[pagewise]{lineno}
\usepackage[square,numbers]{natbib}

\usepackage[colorlinks=true, linkcolor=blue, citecolor=blue, urlcolor=blue, unicode=true]{hyperref}
\usepackage{cleveref}

\newtheorem{theorem}{Theorem}[section]

\newtheorem{proposition}[theorem]{Proposition}
\newtheorem{corollary}[theorem]{Corollary}
\newtheorem{conjecture}[theorem]{Conjecture}

\theoremstyle{definition}
\newtheorem{definition}[theorem]{Definition}
\newtheorem{example}[theorem]{Example}
\theoremstyle{remark}
\newtheorem{remark}[theorem]{Remark}
\newtheorem*{remarks}{Remark}

\DeclareMathOperator{\D}{D}
\DeclareMathOperator{\Hessian}{H} 
\DeclareMathOperator{\Tr}{Tr}

\DeclareMathOperator{\diag}{diag}

\DeclareMathOperator{\vecop}{vec}

\newcommand\inp[2]{\left\langle #1,\,#2 \right\rangle}
\newcommand{\norm}[1]{\left\lVert#1\right\rVert}
\newcommand{\abs}[1]{\lvert#1\rvert}
\newcommand{\id}{\mathrm{I}}
\newcommand*{\intzerotoinfty}[2]{\int_0^{+\infty}#1\mathop{}\!\mathrm{d}#2}

\DeclarePairedDelimiterX{\infdivx}[2]{(}{)}{%
  #1\;\delimsize\|\;#2%
}

\title[Long-Step Path-Following Algorithm]{Long-Step Path-Following Algorithm for Quantum Information Theory: Some Numerical Aspects and Applications}
\author[Leonid Faybusovich and Cunlu Zhou]{}
\email{lfaybuso@nd.edu}
\email{czhou3@nd.edu}

\begin{document}

\maketitle

\centerline{\scshape Leonid Faybusovich and Cunlu Zhou}
\medskip
{\footnotesize
 \centerline{Department of Mathematics}
   \centerline{University of Notre Dame}
   \centerline{Notre Dame, IN 46556, USA}
} 


\medskip

\begin{abstract}
We consider some important computational aspects of the long-step path-following algorithm developed in our previous work and show that a broad class of complicated optimization problems arising in quantum information theory can be solved using this approach. In particular, we consider one difficult and important optimization problem in quantum key distribution and show that our method can solve problems of this type much faster in comparison with (very few) available options.
\end{abstract}

\section{Introduction}
In \cite{fayzhou19longstep} we developed a long-step path-following algorithm to deal with a broad class of symmetric programming problems with nonlinear objective functions. In our recent work \cite{fayzhou19entanglement} we noticed that this class includes a number of difficult and important convex optimization problems arising in quantum information theory. One goal of this paper is to present some important computational aspects of our long-step path-following algorithm that are essential for solving optimization problems of this type. In particular, we show how to derive the analytic expressions of several complicated Hessians and their vectorized forms which are particularly important for practical implementation. Unlike \cite{fayzhou19longstep, fayzhou19entanglement}, we completely avoid using the language of Euclidean Jordan algebras and concentrate on examples with semidefinite constraints. For simplifying our argument and the purpose of illustration, we restrict our discussions and numerical experiments on real symmetric matrices throughout the paper, but all the results can be extended to Hermitian matrices within the Jordan algebraic scheme. By doing so, we hope that our results will reach a broader research community. 

In \cref{sec:longstepalg} we first describe the major features of our algorithm in a broad setting and then discuss in detail the structure of the algorithm for two concrete cases involving semidefinite constraints. In particular, calculations of the Newton directions are discussed. In \cref{sec:matrixmonotone} we consider a class of nonlinear objective functions constructed with the help of matrix monotone functions. Several classes of optimization problems arising in quantum information theory fit into this class, including the \emph{relative entropy of entanglement}, the \emph{fidelity function} used in quantum tomography, and the \emph{relative R\'{e}nyi entropy function}. The structure of the arising Hessians is described in detail, and several important practical aspects for implementation are discussed including vectorization. In \cref{sec:qkd} we consider another important and particularly difficult optimization problem in quantum key distribution (QKD) \cite{norbertqkd, norbertqkd2}. The objective function involves the so-called \emph{quantum relative entropy} function. To the best of our knowledge, there are so far no efficient algorithms available for this type of problem. Generic first-order methods (e.g., the Frank–Wolfe algorithm) are used in \cite{norbertqkd,norbertqkd2}, but the convergence in general is slow and unstable. A more robust method developed in \cite{fawzi2018cvxquad,fawziparrilo18} can be applied, but due to its inherent complexity the method quickly becomes unusable as the problem size increases as shown in \cref{table:qkd}. Although at this stage we have yet been able to establish the compatibility condition \eqref{eq:compcondcase2} for the quantum relative entropy function, we demonstrate that in principle our long-step path-following algorithm can be used to solve the QKD problem efficiently, and the numerical results are indeed quite stunning. For example, compared to the approach of \cite{fawzi2018cvxquad, fawziparrilo18}, their method simply cannot solve the problem for most of the testing examples, and for one of the cases it solves, our method is $9000$ times faster! We conclude our paper in \cref{sec:conclusion} by discussing some future work regarding possible generalizations and improvements of our approach.

\section{Long-Step Path-Following Algorithm}\label{sec:longstepalg}
Let ($\mathbb{E},\,\inp{}{}$) be a Euclidean real vector space with scalar product $\inp{}{}$. Let $a + \mathcal{X}$ be an affine space of $\mathbb{E}$, $\Omega \subseteq \mathbb{E}$ an open convex set and $\overline{\Omega}$ its closure. Let $B(x)$ be the standard self-concordant barrier on $\Omega$ with barrier parameter $r$:
\begin{equation}
\begin{split}
        \abs{\D^3 B(x)(\xi,\xi,\xi)} \, &\leq \, 2\left[\D^2B(x)(\xi,\xi)\right]^{\frac{3}{2}},\\
    \sup_{\xi \in \mathbb{E}}[2\inp{\nabla B(x)}{\xi} &- \inp{\Hessian_B(x)(\xi)}{\xi}] \leq \eta,\\
     x & \in\Omega, \; \xi\in \mathbb{E},\\
     r &= \eta_{\min}.\\
\end{split}
\end{equation}
For a comprehensive discussion of self-concordant barriers, see \cite[section~5.3]{nesterovNewBook}.

Recall that the gradient $\nabla f(x) \in \mathbb{E}$ and Hessian $\Hessian_f(x):\mathbb{E} \to \mathbb{E}$ are defined as follows:
    \begin{equation}
    \begin{split}
        \D f(x)(\xi) &= \inp{\nabla f(x)}{\xi}, \ x\in\Omega,\;\xi\in \mathbb{E},\\
        \D^2 f(x)(\xi,\eta) &= \inp{\Hessian_f(x)\xi}{\eta}, \ x\in\Omega,\;\xi,\eta\in \mathbb{E},
    \end{split}
    \end{equation}
where $\D^k f(x)$ is the $k$-th Fr\'{e}chet derivative of $f$ at $x$. We denote by $C^k(\Omega)$ the vector space of $k$-times continuously differentiable real-valued functions on $\Omega$.

\begin{definition}\label{def:selfconcordant}
Let $F:\Omega \to \mathbb{R}$, $F \in C^3(\Omega)$, be a convex function on $\Omega$. We say that $F$ is \emph{$\kappa$-self-concordant}, if $\exists\,\kappa \geq 0$ such that 
\begin{equation}\label{eq:1}
    \abs{\D^3F(x)(\xi,\xi,\xi)} \, \leq \, 2\kappa\left[\D^2F(x)(\xi,\xi)\right]^{\frac{3}{2}}, \  x\in\Omega, \; \xi\in \mathbb{E}. 
\end{equation}
\end{definition}
We assume that 
\begin{equation}\label{eq:2}
    F(x) \to +\infty,\ x \to \partial \Omega.
\end{equation}
We also assume that the Hessian $\Hessian_F(x)$ is a positive definite symmetric linear operator on $\mathbb{E}$ for all $x\in\Omega$. Given $\xi\in \mathbb{E}$, $x\in\Omega$, we define a norm 
    \begin{equation}\label{eq:3}
        \norm{\xi}_x = \inp{\Hessian_F(x)\xi}{\xi}^{\frac{1}{2}} = [\D^2F(x)(\xi,\xi)]^{\frac{1}{2}}.
    \end{equation}
Under assumptions of \cref{def:selfconcordant} and \eqref{eq:2}, at any point $x \in \Omega$, there exists a so-called \emph{Dikin ellipsoid} inside $\Omega$ \cite[theorem~5.1.5]{nesterovNewBook}:
    \begin{equation}\label{eq:dikin}
        W_s(x) = \{y\in\Omega : \norm{y-x}_x \leq s\} \subset \Omega,\ \forall s<\frac{1}{\kappa}.
    \end{equation}

Now consider the following convex programming problem:
\begin{equation}\label{eq:conicprog}
    \begin{aligned}
    f(x) &\to \min,\\
    x \in \overline{\Omega} &\cap (a + \mathcal{X}),
    \end{aligned}
\end{equation}
where $f \in C^3(\Omega)$, $f$ is continuous on $\overline{\Omega}$ and convex. We assume that the feasible set is bounded and has a nonempty (relative) interior. 

\begin{definition}\label{def:comp}
 $f$ is said to be  \emph{$\nu$-compatible} with $B(x)$ if $\exists\,\nu \geq 1$ such that
\begin{equation}\label{eq:nucomp}
    \abs{\D^3f(x)(\xi,\xi,\xi)} \leq \nu \D^2f(x)(\xi,\xi)[\D^2B(x)(\xi,\xi)]^{\frac{1}{2}},\ \forall\xi\in \mathbb{E}. 
\end{equation}
\end{definition}

Subsequent results are proved\footnote{More precisely, the results are proved for the special case when $\kappa = 1$ and $\nu = 2$ in \cite{fayzhou19longstep}, to which the general case can be reduced by the standard procedure of rescaling the barrier function (see, e.g., \cite[p.~367]{nesterovNewBook}).} in \cite{fayzhou19longstep} for $\Omega$ being the cone of squares of a Euclidean Jordan algebra and $B(x),\,x\in\Omega$, being the standard barrier. However, for the understanding of this paper, no knowledge of Jordan algebras is assumed. We will formulate two special cases in \cref{sec:specialcase1,sec:specialcase2} and show concrete calculations of the so-called \emph{Newton direction} and \emph{Newton decrement}. 

\begin{proposition}\label{prop:comp}
    Let $f$ be $\nu$-compatible with $B(x)$. Then 
    \begin{equation}
        F_\beta(x) = \beta f(x) + B(x),\ x\in\Omega,\,\beta\geq 0,
    \end{equation}
    is $\left(1+\frac{\nu}{3}\right)$-self-concordant, i.e., 
    \begin{equation}
    \abs{\D^3F_{\beta}(x)(\xi,\xi,\xi)} \, \leq \, 2\left(1+\frac{\nu}{3}\right)\left[\D^2F_{\beta}(x)(\xi,\xi)\right]^{\frac{3}{2}}, \  x\in\Omega, \; \xi\in \mathbb{E}. 
\end{equation}
\end{proposition}
This is a straightforward extension of Lemma A.2 in \cite[Appendix~A]{hertogbook} to the general standard barrier $B(x)$, and we omit the proof here. 

With our early assumptions of $f$ and the feasible set, it is easy to see that $F_{\beta}(x)$ has a unique minimum $x(\beta)$ for each $\beta \geq 0$. Now we introduce the Newton direction $p_\beta(x)$ of $F_\beta$ at $x\in(a + \mathcal{X})\cap\Omega$:
    \begin{equation}\label{eq:newtondir}
    \begin{split}
        \Hessian_{F_\beta}(x)p_\beta(x) &= - (\nabla F_{\beta}(x)+\mu_\beta(x)),\\
        \mu_\beta(x) &\in \mathcal{X}^{\perp},\ p_\beta(x) \in \mathcal{X},
    \end{split}
    \end{equation}
    and the Newton decrement of $F_\beta$ at $x$:
        \begin{equation}\label{eq:41}
            \delta_\beta(x) \triangleq \inp{p_\beta(x)}{\Hessian_{F_\beta}(x)p_\beta(x)}^{\frac{1}{2}},\ x\in (a + \mathcal{X})\cap\Omega.
        \end{equation}
    Note that
    \begin{equation}\label{eq:42}
        \delta_\beta(x)^2 = -\inp{\nabla F_\beta(x)}{p_\beta(x)}.
    \end{equation}
Under the assumption of \cref{prop:comp}, we have the following results. 
\begin{proposition}
Given $x\in(a + \mathcal{X})\cap\Omega$, let $\delta_\beta(x)\leq \frac{1}{3\kappa}$. Then
    \begin{equation}
        F_\beta(x) - F_\beta(x(\beta)) \leq \frac{\delta_\beta(x)^2}{1-[\frac{9}{4}\kappa\delta_\beta(x)]^2}.
    \end{equation}
\end{proposition}

\begin{proposition}
    Let $r$ be the barrier parameter of $B(x)$. Given $x\in(a + \mathcal{X})\cap\Omega$ and $\delta_\beta(x) \leq \frac{1}{3\kappa}$, we have
    \begin{equation}
        \abs{f(x)-f(x(\beta))} \leq \left[\frac{\delta_\beta(x)}{1-\frac{9}{4}\kappa\delta_\beta(x)}\cdot\frac{1+\kappa\delta_\beta(x)^2}{1 - \kappa\delta_\beta(x)}\right]\frac{\sqrt{r}}{\beta}.
    \end{equation}
\end{proposition}

\begin{algorithm}[H]
\DontPrintSemicolon
Set $\beta_0 > 0$, and $\theta > 0$. Choose an accuracy $\epsilon > 0$ and find an initial point $x_0 \in (a + \mathcal{X})\cap\Omega$ such that\footnote{For example, use the so-called \emph{analytic center} \cite[Definition~5.3.3]{nesterovNewBook}.}
$$\delta_{\beta_0}(x_0) \leq \frac{1}{3\kappa}.$$\;
\vspace{-5mm}
At $i$-th (outer) iteration ($i \geq 0$), set 
$$\beta_{i+1} = (1+\theta)^{i+1}\beta_0.$$
Find $x_{i+1} \in (a + \mathcal{X})\cap\Omega$ such that $\delta_{\beta_{i+1}}(x_{i+1})\leq \frac{1}{3\kappa}$ by performing several Newton steps (inner iteration) for the function $F_{\beta_{i+1}}$, using $x_{i}$ as the starting point:
$$x_{i} = x_{i} + \alpha p_{\beta_{i+1}}(x_{i}).$$\;
\vspace{-5mm}
Stop the process if 
$$\beta_i \geq \frac{4r}{\epsilon}.$$\;
\vspace{-5mm}
 \caption{A Long-Step Path-Following Algorithm}
\end{algorithm}
\begin{remark}
Note that in the second step, $\alpha$ is obtained by performing a line search (see, e.g., \cite{hertog92classical}) for each inner iteration, where $0 < \alpha < \alpha_{\max}$, and $\alpha_{\max}$ is the largest positive number such that 
$x_{i} + \alpha p_{\beta_{i+1}}(x_{i})$ stays feasible\footnote{Note that for Newton direction $p_{\beta_{i+1}}(x_{i})$ we always have $x_{i} + \alpha p_{\beta_{i+1}}(x_{i})\in a + \mathcal{X}$, and the feasibility is mainly about membership to $\Omega$.}. We used a simple binary-search with gradient in our implementation. 
\end{remark}

\subsection{Complexity Estimates}

\begin{theorem}\label{thm:1}
    Given $\epsilon > 0$ and
    \begin{equation}
        i\geq\frac{\ln(\frac{4r}{\epsilon\beta_0})}{\ln(1+\theta)}.
    \end{equation}
    Then
    $$f(x_{i})-f(x^*)\leq\epsilon,$$
    where $x^*$ is an optimal solution to the problem \eqref{eq:conicprog}. 
\end{theorem}

\begin{theorem}\label{thm:2}
    Each outer iteration requires at most 
    \[
        \frac{22}{3} + 22\theta\left(\frac{5}{2}\kappa\sqrt{r} + \frac{\theta\kappa^2 r}{\theta + 1}\right)
    \]
    inner iterations.
\end{theorem}

Taking into account \cref{thm:1,thm:2}, we get the following complexity result for our algorithm. 
\begin{theorem}\label{thm:3}
    An upper bound for the total number of Newton iterations is given by
    \[
        \frac{\ln(\frac{4r}{\epsilon\beta_0})}{\ln(1+\theta)}\left(\frac{22}{3} + 22\theta\left(\frac{5}{2}\kappa\sqrt{r} + \frac{\theta\kappa^2 r}{\theta + 1}\right)\right).
    \]
\end{theorem}

\subsection{Two Special Cases}\label{sec:twospecialcases}
We consider two special cases of the convex programming problem \eqref{eq:conicprog}. Without loss of generality, let $\mathbb{E} = \mathbb{S}^n$, the real vector space of $n \times n$ symmetric matrices. We denote by $\mathbb{S}^n_{+}$ and $\mathbb{S}^n_{++}$ the convex cone of positive semidefinite matrices and positive definite matrices respectively. Let $\mathbb{R}^N_+$ denote the nonnegative orthant of $\mathbb{R}^N$. 

We use $A \succeq 0$ for positive semidefiniteness (i.e., $x^T A x \geq 0,\,\forall x \in \mathbb{R}^n$) and $A \succ 0$ for positive definiteness (i.e., $x^T Ax > 0,\ \forall x \in \mathbb{R}^n\setminus \{0\}$). Furthermore, we use notations
$$A \succeq B\text{ if }A-B \succeq 0,$$
and
$$A \succ B \text{ if } A-B \succ 0.$$
The scalar product $\inp{A}{B}$, $A,B \in \mathbb{S}^n$, is defined as $$\inp{A}{B} = \Tr(A^TB) = \sum_{i,j}A_{ij}B_{ij}.$$ 

Recall the following facts about $B(X) = -\ln\det(X),\,X\in\mathbb{S}^n_{++}$:
\begin{align}
    \nabla B(X) &= - X^{-1},\\
    \Hessian_B(X) &= P(X^{-1}),
\end{align}
where 
\begin{equation}\label{eq:pop}
    P(X^{-1})Y = X^{-1}YX^{-1},\,\forall Y\in\mathbb{S}^n.
\end{equation}

\subsection{Type I}\label{sec:specialcase1} We consider the following optimization problem:
\begin{equation}\label{eq:case1}
    \begin{aligned}
    f(X) &\to \min,\\
    \inp{A_i}{X} &\leq b_i,\,i=1,\ldots,m,\\
    \inp{A_i}{X} &= b_i,\,i=m+1,\ldots,N,\\
    X &\succeq 0,
    \end{aligned}
\end{equation}
where we assume $A_i$'s are linearly independent. We also assume the feasible set is bounded and has a nonempty (relative) interior. 

Note that we can rewrite \eqref{eq:case1} in the following form:
\begin{equation}\label{eq:case1prime}
    \begin{aligned}
    f(X) &\to \min,\\
    \inp{A_i}{X} + x_i &= b_i,\,i=1,\ldots,m,\\
    \inp{A_i}{X} &= b_i,\,i=m+1,\ldots,N,\\
    x_i &\geq 0,\,i=1,\ldots,m,\\
    X &\succeq 0.
    \end{aligned}
\end{equation}
Notice that the closed convex cone $\overline{\Omega}$ in \eqref{eq:conicprog} is $\mathbb{S}^n_{+} \times \mathbb{R}^{m}_{+}$ in this case.

Recall \cref{def:comp} and \cref{prop:comp}. The $\nu$-compatibility ($\nu \geq 1$) condition for $f(x)$ in \eqref{eq:case1prime} and the corresponding auxiliary barrier family of optimization problems are given as follows: 
\begin{equation}\label{eq:compcondcase1}
    \abs{\D^3f(X)(\xi,\xi,\xi)} \leq \nu \D^2f(X)(\xi,\xi)[\D^2 B(X)(\xi,\xi)]^{\frac{1}{2}},\ \forall\xi\in \mathbb{S}^n,
\end{equation}
and  
\begin{equation}\label{eq:auxcase1}
    \begin{aligned}
    F_{\beta}(X;x) &= \beta f(X) -\ln\det(X) - \sum_{i=1}^m \ln(x_i) \to \min,\\
    \inp{A_i}{X} + x_i &= b_i,\,i=1,\ldots,m,\\
    \inp{A_i}{X} &= b_i,\,i=m+1,\ldots,N,\\
    x_i &> 0,\,i=1,\ldots,m,\\
    X &\succ 0.
    \end{aligned}
\end{equation}

Next we show how to obtain the Newton direction and the Newton decrement. First we can rewrite the constraints into the following compact form\footnote{Strictly speaking here we should use $\vecop(A_i)$ and $\vecop(X)$, where the $\vecop(\cdot)$ operator is only introduced at a later time. Similar note for the gradient and the Newton direction discussed below.}:

\begin{equation}
    \inp{\begin{bmatrix}A_i\\e_i\end{bmatrix}}{\begin{bmatrix}X\\x\end{bmatrix}} = \inp{A_i}{X} + \inp{e_i}{x} = b_i,\,i=1,\ldots,m,
\end{equation}
and 
\begin{equation}
    \inp{A_i}{X} = b_i,\,i=m+1,\ldots,N,
\end{equation}
where $e_i = [0,\ldots,1,\ldots,0]^T$ (with $1$ at the $i$-th position).

We have the gradient 
\begin{equation}\label{eq:case1grad}
    \nabla F_{\beta}(X;x) = \begin{bmatrix}\nabla F_{\beta}^{(1)}(X;x)\\\nabla F_{\beta}^{(2)}(X;x)\end{bmatrix},
\end{equation}
where 
\begin{equation*}
    \nabla F_{\beta}^{(1)}(X;x) = \beta \nabla f(X) - X^{-1},
\end{equation*}
and 
\begin{equation*}
    \nabla F_{\beta}^{(2)}(X;x) = \left[-\frac{1}{x_1},\ldots,-\frac{1}{x_m}\right]^T.
\end{equation*}
We have the Hessian
\begin{equation}\label{eq:case1hessian}
    \Hessian_{F_{\beta}}(X;x) = \begin{bmatrix}\Hessian^{(1)}_{F_{\beta}}(X;x) & \mathbf{0}\\\mathbf{0} & \Hessian^{(2)}_{F_{\beta}}(X;x)\end{bmatrix},
\end{equation}
where
\begin{equation*}
    \Hessian^{(1)}_{F_{\beta}}(X;x) = \beta \Hessian_f(X) + P(X^{-1}),
\end{equation*}
and
\begin{equation*}
    \Hessian^{(2)}_{F_{\beta}}(X;x) = \diag\left( \left[\frac{1}{x^2_1},\ldots,\frac{1}{x^2_m}\right]\right).
\end{equation*}
From \eqref{eq:newtondir}, we have 
\begin{equation}\label{eq:case1newtondir1}
    \Hessian_{F_{\beta}}(X;x)p_{\beta}(X;x) = - \nabla F_{\beta}(X;x) + \sum_{j=1}^m \lambda_j \begin{bmatrix}A_j\\e_j\end{bmatrix} + \sum_{j=m+1}^N\lambda_j \begin{bmatrix}A_j\\0\end{bmatrix},
\end{equation}
where $p_{\beta}(X;x)$ is the Newton direction:
\begin{equation*}
    p_{\beta}(X;x) = \begin{bmatrix}p^{(1)}_{\beta}(X;x)\\p^{(2)}_{\beta}(X;x)\end{bmatrix},
\end{equation*}
which satisfies the following conditions:
\begin{equation}\label{eq:case1newtondir2}
    \inp{\begin{bmatrix}A_i\\e_i\end{bmatrix}}{\begin{bmatrix}p^{(1)}_{\beta}(X;x)\\p^{(2)}_{\beta}(X;x)\end{bmatrix}} = 0,\, i = 1,\ldots,m,
\end{equation}
and
\begin{equation}\label{eq:case1newtondir3}
    \inp{A_i}{p^{(1)}_{\beta}(X;x)} = 0,\, i = m+1,\ldots,N.
\end{equation}

From \eqref{eq:case1newtondir1}, we have 
\begin{equation}\label{eq:case1newdir1}
    p^{(1)}_{\beta}(X;x) = \Hessian^{(1)}_{F_{\beta}}(X;x)^{-1} \left[-\nabla F_{\beta}^{(1)}(X;x) + \sum_{j=1}^N \lambda_j A_j\right],
\end{equation}
and
\begin{equation}\label{eq:case1newdir2}
    \begin{split}
        p^{(2)}_{\beta}(X;x) &= \Hessian^{(2)}_{F_{\beta}}(X;x)^{-1} \left[-\nabla F_{\beta}^{(2)}(X;x) + \sum_{j=1}^m \lambda_j e_j\right]\\
        &= \begin{bmatrix}x_1 + \lambda_1 x^2_1\\\vdots\\x_m + \lambda_m x^2_m\end{bmatrix}.
    \end{split}
\end{equation}

Combining \eqref{eq:case1newtondir2}, \eqref{eq:case1newtondir3}, \eqref{eq:case1newdir1}, and \eqref{eq:case1newdir2}, we can assemble a linear system of equations:

\begin{multline}\label{eq:case1eq1}
    \sum_{j=1}^N \lambda_j \inp{A_i}{\Hessian^{(1)}_{F_{\beta}}(X;x)^{-1}A_j} + \lambda_i x_i^2 \\= \inp{A_i}{\Hessian^{(1)}_{F_{\beta}}(X;x)^{-1}\nabla F_{\beta}^{(1)}(X;x)} - x_i,\,i=1,\ldots,m,
\end{multline}
and
\begin{multline}\label{eq:case1eq2}
    \sum_{j=1}^N \lambda_j \inp{A_i}{\Hessian^{(1)}_{F_{\beta}}(X;x)^{-1}A_j} = \inp{A_i}{\Hessian^{(1)}_{F_{\beta}}(X;x)^{-1}\nabla F_{\beta}^{(1)}(X;x)},\,i= m+1,\ldots,N.
\end{multline}

By solving \eqref{eq:case1eq1} and \eqref{eq:case1eq2}, we obtain $\lambda_j,\,j=1,\ldots,N$, and hence the Newton direction $p_{\beta}(X;x)$. By \eqref{eq:42}, we get the Newton decrement
\begin{equation}
    \begin{split}
        \delta_{\beta}(X;x) &= \sqrt{-\inp{p_{\beta}(X;x)}{\nabla F_{\beta}(X;x)}}\\
        &= \sqrt{-\left[\inp{p^{(1)}_{\beta}(X;x)}{\nabla F^{(1)}_{\beta}(X;x)} + \inp{p^{(2)}_{\beta}(X;x)}{\nabla F^{(2)}_{\beta}(X;x)}\right]}.\\
    \end{split}
\end{equation}

\subsection{Type II}\label{sec:specialcase2} We consider optimization problems of the following form:
\begin{equation}\label{eq:case2}
    \begin{aligned}
    f(X) &+ g(Y) \to \min,\\
    \inp{A_i}{X} &= b_i,\,i=1,\ldots,m,\\
    Y &= \mathcal{L}(X),\\
    X &\succeq 0,\\
    Y &\succeq 0,
    \end{aligned}
\end{equation}
where $\mathcal{L}: \mathbb{S}^n_+ \to \mathbb{S}^k_+$ is some linear operator. We assume the feasible set is bounded and has a nonempty (relative) interior. We further assume that $A_i$'s are linearly independent. Note that $\overline{\Omega} = \mathbb{S}^n_{+} \times \mathbb{S}^k_{+}$ in this case.

The $\nu$-compatibility ($\nu \geq 1$) condition for \eqref{eq:case2} and the corresponding auxiliary barrier family of optimization problems are given as follows: 
\begin{multline}\label{eq:compcondcase2}
\abs{\D^3f(X)(\xi,\xi,\xi) + \D^3 g(Y)(h,h,h)} \\
\leq \nu \left[\D^2 f(X)(\xi,\xi) + \D^2 g(Y)(h,h)\right] \left[\D^2 B_1(X)(\xi,\xi) + \D^2 B_2(Y)(h,h)\right]^{\frac{1}{2}},\,
\forall\xi, h \in \mathbb{S}^n,
\end{multline}
where $B_1(X) = -\ln\det(X)$ and $B_2(Y) = -\ln\det(Y)$,
and 
\begin{equation}\label{eq:auxcase2}
    \begin{aligned}
    F_{\beta}(X) &= \beta (f(X) + g(Y)) -\ln\det(X) - \ln\det(Y) \to \min,\\
    \inp{A_i}{X} &= b_i,\,i=1,\ldots,m,\\
    Y &= \mathcal{L}(X),\\
    X &\succ 0,\\
    Y &\succ 0. 
    \end{aligned}
\end{equation}

Equivalently, we can rewrite \eqref{eq:auxcase2} as follows
\begin{equation}\label{eq:auxcase2prime}
    \begin{aligned}
    F_{\beta}(X) &= \beta (f(X) + g(\mathcal{L}(X))) -\ln\det(X) - \ln\det(\mathcal{L}(X)) \to \min,\\
    \inp{A_i}{X} &= b_i,\,i=1,\ldots,m,\\
    X &\succ 0,\\
    \mathcal{L}(X) &\succ 0. 
    \end{aligned}
\end{equation}

Next we show calculations of the Newton direction $p_{\beta}(X)$ and the Newton decrement $\delta_{\beta}(X)$. By \eqref{eq:newtondir}, we have
\begin{align}
    \Hessian_{F_{\beta}}(X)p_{\beta}(X) &= - \nabla F_{\beta}(X) + \sum_{j=1}^m \lambda_j A_j,\\
    \inp{A_i}{p_{\beta}(X)} &= 0,\,i=1,\ldots,m,
\end{align}
by which we can assemble the following linear system of equations:
\begin{equation}\label{eq:case2linearsys}
    \sum_{j=1}^m \lambda_j \inp{A_i}{\Hessian_{F_{\beta}}(X)^{-1}A_j} = \inp{A_i}{\Hessian_{F_{\beta}}(X)^{-1}\nabla F_{\beta}(X)},\,i=1,\ldots,m.
\end{equation}
By solving \eqref{eq:case2linearsys} we obtain $\lambda_j,\,j=1,\ldots,m$, and hence the Newton direction
\begin{equation}\label{eq:case2newtondir}
    p_{\beta}(X) = \Hessian_{F_{\beta}}(X)^{-1} \left(- \nabla F_{\beta}(X) + \sum_{j=1}^m \lambda_j A_j\right),
\end{equation}
and the Newton decrement
\begin{equation}
    \delta_{\beta}(X) = \sqrt{-\inp{\nabla F_{\beta}(X)}{p_{\beta}(X)}}.
\end{equation}

For calculations of $\nabla F_{\beta}(X)$ and $\Hessian_{F_{\beta}}(X)$, let $h(X) = g(\mathcal{L}(X))$ and $\zeta(X) = B(\mathcal{L}(X)) = - \ln\det(\mathcal{L}(X))$. 
We have (by chain rule) for all $\xi \in \mathbb{S}^n$,
\begin{equation}
    \begin{split}
        \D h(X)(\xi) &= \D g(\mathcal{L}(X))(\mathcal{L}(\xi))\\
        &= \inp{\nabla g(\mathcal{L}(X))}{\mathcal{L}(\xi)}\\
        &= \inp{\mathcal{L}^T\nabla g(\mathcal{L}(X))}{\xi}\\
        &= \inp{\nabla h(X)}{\xi},
    \end{split}
\end{equation}
which implies that
\begin{equation}\label{eq:gradgeneral}
    \nabla h(X) = \mathcal{L}^T\nabla g(\mathcal{L}(X)).
\end{equation}
We further have
\begin{equation}\label{eq:hessiangeneral1}
    \begin{split}
        \D^2 h(X)(\xi,\xi) &= \D^2 g(\mathcal{L}(X))(\mathcal{L}(\xi),\mathcal{L}(\xi))\\
        &= \inp{\Hessian_g(\mathcal{L}(X))\mathcal{L}(\xi)}{\mathcal{L}(\xi)}\\
        &= \inp{\mathcal{L}^T\Hessian_g(\mathcal{L}(X))\mathcal{L}(\xi)}{\xi}\\
        &= \inp{\Hessian_h(X)\xi}{\xi},
    \end{split}
\end{equation}
which implies that
\begin{equation}\label{eq:hessiangeneral}
    \Hessian_h(X)(\xi) = \mathcal{L}^T\Hessian_g(\mathcal{L}(X))\mathcal{L}(\xi).
\end{equation}
Similarly, we have
\begin{equation}
    \nabla \zeta(X) = \mathcal{L}^T\nabla B(\mathcal{L}(X)) = \mathcal{L}^T(-\mathcal{L}(X)^{-1}). 
\end{equation}
and
\begin{equation}
    \Hessian_{\zeta}(X) = \mathcal{L}^T\Hessian_B(\mathcal{L}(X))\mathcal{L} = \mathcal{L}^T P(\mathcal{L}(X)^{-1})\mathcal{L}, 
\end{equation}
where $P(\cdot)$ was introduced in \eqref{eq:pop}.

\section{Matrix Monotone Functions}\label{sec:matrixmonotone}
Let 
$$g:[0,+\infty) \to \mathbb{R}$$
 be a real-valued function. We say that $g$ is \emph{matrix monotone} (\emph{anti-monotone}) if for any real symmetric matrices of the same size such that $A \succeq 0$, $B \succeq 0$ and $A \succeq B$, we have
$$g(A) \succeq g(B)\ (g(A) \preceq g(B)).$$
It is obvious that if $g$ is matrix monotone then $-g$ is matrix anti-monotone and vice versa. 
In \cite{fayzhou19entanglement}, we proved that for any matrix anti-monotone function 
$$g:[0,+\infty) \to \mathbb{R},$$ 
we have the following compatibility result (adapted for the case of symmetric matrices).
\begin{theorem}\label{thm:monotonecomp}
Let $C \in \mathbb{S}^n_{+}$ and $B(X) = - \ln\det(X),\,X\in\mathbb{S}^n_{++}$. Then
\begin{equation*}
    \abs{\D^3 \varphi_c(X)(\xi,\xi,\xi)} \leq 3 \D^2 \varphi_c(X)(\xi,\xi) \sqrt{\D^2 B(X)(\xi,\xi)},\ X\in \mathbb{S}^n_{++},\,\xi\in \mathbb{S}^n,
\end{equation*}
where 
\begin{equation}
    \varphi_c(X) = \inp{C}{g(X)}.
\end{equation}
\end{theorem}
Hence, by \cref{prop:comp}, we have the following self-concordance result.
\begin{corollary}\label{cor:monotonecon}
For any $\beta \geq 0$, the function
$$\Phi_{\beta}(X) = \beta \varphi_c(X) + B(X),\ X\in\mathbb{S}^n_{++},$$
is $\kappa$-self-concordant on $\mathbb{S}^n_{++}$ with $\kappa = 2$.
\end{corollary}
With \cref{thm:monotonecomp} and \cref{cor:monotonecon}, the long-step path-following algorithm discussed in \cref{sec:longstepalg} can then be applied to optimization problems involving objective functions of the form:
\begin{equation}\label{eq:objfun}
    \varphi_c(X) = \inp{C}{g(X)} = \Tr(Cg(X)),\,C \succeq 0.
\end{equation}

For implementation, we show how the analytic expressions of the gradient and the Hessian of $\varphi_c(X)$ can be obtained. Let
$$X = U\Lambda U^T$$
be a spectral decomposition of $X$, where $\Lambda = \diag(\lambda_1,\ldots,\lambda_n)$ and $UU^T = \id$. 

For a continuously differentiable function $h:[0,+\infty) \to \mathbb{R}$, we introduce the first divided difference $h^{[1]}$:
\begin{equation}\label{eq:1stdiv}
 h^{[1]}(\lambda_i,\lambda_j)
    = \left\{\begin{aligned}
             \frac{h(\lambda_i) - h(\lambda_j)}{\lambda_i - \lambda_j}&, &  &\lambda_i \ne \lambda_j, \\
             h'(\lambda_i)&, &  &\lambda_i = \lambda_j,
        \end{aligned}\right.
\end{equation}
and the second divided difference $h^{[2]}$:
\begin{equation}\label{eq:2nddd}
    h^{[2]}(\lambda_i, \lambda_j, \lambda_k) = \frac{h^{[1]}(\lambda_i,\lambda_j) - h^{[1]}(\lambda_i,\lambda_k)}{\lambda_j - \lambda_k}
\end{equation}
for distinct $\lambda_i$, $\lambda_j$, and $\lambda_k$,  while for other cases the function is defined by taking limits in \eqref{eq:2nddd}, e.g., 
\begin{equation*}
    h^{[2]}(\lambda, \lambda, \lambda) = \frac{1}{2}h''(\lambda).
\end{equation*}
Lastly, recall the Schur product for $m \times n$ matrices $A$ and $B$ is defined as 
$$[A \circ B]_{ij} = A_{ij}B_{ij},$$
and the vectorization operator $\vecop(\cdot)$ for an $n\times m$ matrix $A = [a_{ij}]$:
$$\vecop(A) = [a_{11},\ldots,a_{n1},a_{12},\ldots,a_{n2},\ldots,a_{1m},\ldots,a_{nm}]^T.$$
The following identity is particularly useful\footnote{For complex matrices, we have $\vecop(XYZ^*) = (X\otimes \overline{Z})\vecop(Y)$, where $\overline{Z}$ is the conjugate matrix of $Z$.}:
\begin{equation}\label{eq:vecid}
    \vecop(XYZ^T) = (X\otimes Z)\vecop(Y).
\end{equation}

In \cite{fayzhou19entanglement}, we showed the derivation of the gradient and Hessian of $\varphi_c(X)$ for the case when $g(t)=-\ln(t),\,t>0$. Since the results are derived from the integral representation of $\ln(X)$ and any matrix monotone function admits such an integral representation, we can derive the gradient and Hessian of $\varphi_c(X)$ for any matrix anti-monotone function $g(t),\,t\geq 0$, and the only difference will be  calculations of the first and second divided differences. For completeness, we have reproduced the results for $g(t) = -\ln(t),\,t>0$, in \cref{sec:appendix}.

By \eqref{eq:dtrclnx} we have\footnote{Note the subtle difference that $h(\lambda) = \ln(t),\,t>0,$ in \eqref{eq:dtrclnx} which explains the minus sign there.} 
\begin{equation}
    \D \varphi_c(X)(\xi) = \inp{U\left((U^T C U)\circ g^{[1]}(\Lambda)\right)U^T}{\xi},\ \forall \xi \in \mathbb{S}^n,
\end{equation}
and hence 
\begin{equation}\label{eq:gradmonotone}
    \nabla \varphi_c(X) = U\left((U^T C U)\circ g^{[1]}(\Lambda)\right)U^T,
\end{equation}
where $g^{[1]}(\Lambda)$ is the $n\times n$ first divided difference matrix with $[g^{[1]}(\Lambda)]_{ij} = g^{[1]}(\lambda_i,\lambda_j)$. 

Furthermore, by \eqref{eq:vecid} we have
\begin{equation}\label{eq:gradmonotonevec}
\begin{split}
        \vecop({\nabla\varphi_c(X)}) 
        &= (U\otimes U)\left(\vecop(U^T C U)\circ \vecop(g^{[1]}(\Lambda))\right),\\
        &= (U\otimes U)\left(\diag(\vecop(g^{[1]}(\Lambda)))\vecop(U^T C U)\right)\\
        &= (U\otimes U)\left(\diag(\vecop(g^{[1]}(\Lambda)))(U^T \otimes U^T) \vecop(C)\right)\\
        &= (U \otimes U) \diag(\vecop(g^{[1]}(\Lambda))) (U \otimes U)^T \vecop(C).
\end{split}
\end{equation}

By \eqref{eq:htrclnx} we have
\begin{equation}\label{eq:hessianmonotone}
    \Hessian_{\varphi_c}(X)(\xi) = U \left(\intzerotoinfty{(D\tilde{C}D\tilde{\xi}D + D\tilde{\xi}D\tilde{C}D)}{t}\right) U^T,\ \forall \xi \in \mathbb{S}^n,
\end{equation}
where $D = (\Lambda + t\id)^{-1}$, $\tilde{\xi} = U^T \xi U$, and $\tilde{C} = U^T C U$.

Again by \eqref{eq:vecid} we get
\begin{equation}\label{eq:hessianmonotonevecdev}
    \begin{split}
        \vecop\left(\Hessian_{\varphi_c}(X)(\xi)\right) &= (U\otimes U) \left(\intzerotoinfty{((D\tilde{C}D)\otimes D + D\otimes (D\tilde{C}D))}{t}\right)\vecop(U^T \xi U)\\
        &= (U\otimes U) \left(\intzerotoinfty{((D\tilde{C}D)\otimes D + D\otimes (D\tilde{C}D))}{t}\right) (U\otimes U)^{T} \vecop(\xi),
    \end{split}
\end{equation}
and hence
\begin{equation}\label{eq:hessianmonotonevec}
        \Hessian_{\varphi_c}(X) = (U\otimes U) \left(\intzerotoinfty{((D\tilde{C}D)\otimes D + D\otimes (D\tilde{C}D))}{t}\right) (U\otimes U)^{T}.
\end{equation}

As discussed in \cite{fayzhou19entanglement}, the middle part 
$$S = \intzerotoinfty{((D\tilde{C}D)\otimes D + D\otimes (D\tilde{C}D))}{t}$$
is a sparse block matrix with $(ij,kl)$-th entry: 
\begin{equation}
    S_{ij,kl} = \delta_{kl}\tilde{C}_{ij}\Gamma_{ijl} + \delta_{ij}\tilde{C}_{kl}\Gamma_{jkl},
\end{equation}
where 
$$\delta_{ij} = \left\{\begin{aligned}\,1\quad &\text{ if } i=j,\\ \,0\quad &\text{ if } i \neq j,
\end{aligned}
\right. \text{ and }\Gamma_{ijk} = g^{[2]}(\lambda_i, \lambda_j, \lambda_k),$$
from which we notice that the $ij$-th sub-block matrix is diagonal if $i \neq j$.

Next we illustrate some important numerical aspects of the long-step path-following algorithm through examples. 
\begin{example}
Consider the function 
$$g(t) = t^{-1},\,t>0.$$ 
Clearly, $g$ is matrix anti-monotone. Therefore, we can apply the long-step path-following algorithm to the following optimization problem: for $C \succeq 0$,
\begin{equation}\label{eq:invx}
    \begin{aligned}
    f(X) &= \Tr(CX^{-1})\to \min,\\
    \inp{A_i}{X} &\leq b_i,\,i=1,\ldots,m,\\
    \inp{A_i}{X} &= b_i,\,i=m+1,\ldots,N,\\
    X &\succeq 0.
    \end{aligned}
\end{equation}
Note that this is the optimization problem of type I discussed in \cref{sec:specialcase1}. \Cref{table:invx} shows the numerical results of solving \eqref{eq:invx}. We use the \emph{analytic center} \cite[Definition~5.3.3]{nesterovNewBook} as our initial point, which can be easily obtained by solving the following optimization problem (e.g., use SDPT3 \cite{sdpt3}): 
\begin{equation}\label{eq:miniac}
\begin{aligned}
    f_{ac}(X) &= - \ln\det(X) - \sum_{i=1}^m \ln(x_i) \to \min,\\
    \inp{A_i}{X} + x_i &= b_i,\,i=1,\ldots,m,\\
    \inp{A_i}{X} &= b_i,\,i=m+1,\ldots,N,\\
    X &\succeq 0,\\
    x_i &\geq 0.
\end{aligned}
\end{equation}

Our algorithm is implemented in Matlab, and all the numerical experiments are performed on a personal 15-in Macbook Pro with Intel core i7 and 16 GB memory. Data are randomly generated \footnote{Data can be accessed here: \url{https://doi.org/10.13140/RG.2.2.19100.51847}} and without loss of generality, $\Tr(X) = 1$ is imposed (we assume that the feasible set is bounded). We set $\beta_0 = 0.1$, $\theta = 10$, and $\epsilon = \num{1e-4}$ for all of our tests. In \cref{table:invx}, $nNewton$ is the total number of Newton steps used. $T_{ac}$ is the time for solving the analytic center and $T_{pf}$ is the time for running the long-step path-following algorithm. Both $T_{ac}$ and $T_{pf}$ are averaged over 20 repeated runs. 

\begin{table}[tbhp]
\centering
\footnotesize{
\caption{Numerical Results for \eqref{eq:invx}}\label{table:invx}
\begin{tabular}{ r r r r c r r}
\\\toprule
 & & &\multicolumn{3}{c}{long-step path-following} \\

\midrule
$n$ & $m$ & $N$ & $f_{min}$ & $nNewton$ & $T_{ac}$(s) & $T_{pf}$(s)  \\
\midrule
  4 & 2 & 4 & 27.3538 & 7 & 0.18 & 0.01\\
  8 & 4 & 8 & 8.3264  & 13 & 0.19 & 0.03\\
  16 & 8 & 16 & 18.4274 & 13 & 0.26 & 0.09\\
  32 & 16 & 32 & 39.2516 & 21 & 1.39 & 1.06\\
  64 & 32 & 64 & 91.6534 & 27 & 26.34 & 47.25\\
\bottomrule
\end{tabular}
}
\end{table}

\begin{remarks}
\renewcommand{\labelenumi}{\arabic{enumi}.}
\begin{enumerate}
    \item We noticed that the most time consuming part when running the algorithm is assembling the linear system, e.g., \eqref{eq:case1eq1}, \eqref{eq:case1eq2}, and \eqref{eq:case2linearsys}. 
     \item In general, using vectorization greatly improves  performance and scalability of the algorithm \cite{overton98}.
\end{enumerate}
\end{remarks}
\end{example}

\begin{example}
In quantum information theory, the so-called \emph{relative R\'{e}nyi entropy} is defined as
\begin{equation}
    \begin{aligned}
        \varphi_{\alpha}(X,Y) &= - \Tr(X^{\alpha} Y^{1-\alpha}),\,\alpha \in (0,1),\\
        X,Y &\in \mathbb{S}^n_{++}.
    \end{aligned}
\end{equation}
The function $\varphi_{\alpha}$ is jointly convex in $X,Y$. For a fixed $Y$, the function $X \mapsto  \varphi_{\alpha}(X,Y)$ is matrix anti-monotone, and for a fixed $Y$, the function $X \mapsto  \varphi_{\alpha}(X,Y)$ is matrix anti-monotone. Therefore, for optimization problems involving the relative R\'{e}nyi entropy, our long-step path-following algorithm combined with an alternative minimization procedure (similar to the one used in \cite{fayzhou19entanglement}) can be applied.  
\end{example}

\begin{example}
The relative entropy of entanglement (REE) problem described in \cite{fayzhou19entanglement} involves the following optimization problem which gives a lower bound to the REE of a quantum state $C$ (i.e., $C \succeq 0$ and $\Tr(C) = 1$):
\begin{equation}\label{eq:reeopt}
    \begin{split}
        f(X) &= \Tr(C \ln(C)) - \Tr(C \ln(X)) \to \min,\\
        \Tr(X) &= 1,\\
        \mathcal{L}(X) &\geq 0,\\
        X &\succeq 0,
    \end{split}
\end{equation}
where $\mathcal{L}(\cdot)$ is the so-called \emph{partial transpose} operator. 
Note that the function $\lambda \mapsto -\ln(\lambda)$ is matrix anti-monotone and the REE optimization problem \eqref{eq:reeopt} is of type II discussed in \cref{sec:specialcase2}. Therefore, our long-step path-following is readily to be applied. For numerical results and more details regarding the REE problem, we refer to our previous work \cite{fayzhou19entanglement}.  
\end{example}

\begin{example}
The objective functions based on fidelity \cite{watrous18} have the form
\begin{equation}
    \varphi(X) = - \Tr(\mathcal{L}(X)^{\frac{1}{2}}),
\end{equation}
where $X \in \mathbb{S}^n_{++}$ and $\mathcal{L}(X) = Y^{\frac{1}{2}}XY^{\frac{1}{2}}$ for some fixed $Y \in \mathbb{S}^n_{++}$. Note that the function $\lambda \mapsto -\sqrt{\lambda}$ is matrix anti-monotone. It immediately follows  that our path-following algorithm can be applied to this type of problems as well. 
\end{example}

\subsection{Some Important Observations}
During our numerical experiments, we have the following important observations:
\renewcommand{\labelenumi}{\arabic{enumi}.}
\begin{enumerate}
    \item While conducting extensive numerical experiments with Newton's method (with line search) applied to problems with semidefinite constraints, we noticed a striking difference between the cases of self-concordant and non-self-concordant functions. In the latter case the convergence of Newton's method is rather slow (or there is no convergence to an optimal solution at all) even when the optimal solution lies in the interior of a feasible set. 
    \item For some optimization problems of type II, such as the REE problem \eqref{eq:reeopt}, the barrier term $-\ln\det(X)$ seems unnecessary when running our long-step path-following algorithm. We suspect that this is related to certain properties of the linear operator involved,  and it would be interesting to understand more about this phenomenon. 
\end{enumerate}

\section{An Important Optimization Problem in Quantum Key Distribution}\label{sec:qkd}
Key distribution is used to distribute security keys
to two parties so they can securely share information. While traditional public key distribution is based on the computational intractability of hard mathematical problems, quantum key distribution (QKD) relies on the fundamental law of nature, or more precisely, on the theory of quantum mechanics. QKD has been shown to provide a quantum-secure method of sharing keys which in principle is immune to the power of an eavesdropper \cite{locurty14,scarani09}. 

One of the main theoretical problems in QKD is to calculate the secret key rate for a given QKD protocol, which is essentially to solve the following optimization problem involving the quantum relative entropy function \cite{norbertqkd,norbertqkd2}:
\begin{equation}\label{eq:minirelqkd}
    \begin{aligned}
    f(\tilde{X},\tilde{Y}) &= \Tr(\tilde{X} \ln(\tilde{X})) - \Tr(\tilde{X} \ln(\tilde{Y})) \to \min,\\
    \tilde{X} & = \sum_{j=1}^l K_j X K_j^*,\\
    \tilde{Y} &= \sum_{p=1}^s Z_p \tilde{X} Z_p,\\
    \Tr(A_i X) &= b_i,\,i=1,\ldots,m,\\
    \end{aligned}
\end{equation}
where $X$ is an $n\times n$ density matrix (i.e., $X \succeq 0$ and $\Tr(X) = 1$), $K_j$'s are $k \times n$ matrices such that $\sum_{j=1}^l K_j^* K_j \leq \id$, and $Z_p$'s are $k \times k$ orthogonal projectors such that $\sum_{p=1}^s Z_p = \id$. Note that $k$ usually depends on $n$ (e.g., $k=2n$).

To the best of our knowledge, there are so far no efficient algorithms available for solving \eqref{eq:minirelqkd}. Generic first-order methods (e.g., the Frank–Wolfe algorithm) are used in \cite{norbertqkd,norbertqkd2}, but the convergence is in general slow and unstable. A more robust method developed in \cite{fawzi2018cvxquad,fawziparrilo18} can be applied to \eqref{eq:minirelqkd}, but due to its inherent complexity the method quickly becomes unusable as shown in \cref{table:qkd}. Although at this stage we have yet been able to establish the compatibility condition \eqref{eq:compcondcase2} for the quantum relative entropy function, we will demonstrate that in principle our long-step path-following algorithm can be used to solve \eqref{eq:minirelqkd} efficiently because it is a structured direct method that can exploit all the structural properties of the problem. 

First we can rewrite \eqref{eq:minirelqkd} in the following general form:
\begin{equation}\label{eq:qkdgeneric}
    \begin{aligned}
    f(\tilde{X},\tilde{Y}) &= \Tr(\tilde{X} \ln(\tilde{X})) - \Tr(\tilde{X} \ln(\tilde{Y})) \to \min,\\
    \tilde{X} & = \mathcal{L}_1(X),\\
    \tilde{Y} &= \mathcal{L}_2(X),\\
    \Tr(A_i X) &= b_i,\,i=1,\ldots,m,\\
    X &\succeq 0,
    \end{aligned}
\end{equation}
where $\mathcal{L}_1$ and $\mathcal{L}_2$ are two linear operators of the same type:
\begin{align*}
    \mathcal{L}_1&: X \mapsto \sum_{j=1}^{r_1} K_j X K_j^*,\\
    \mathcal{L}_2&: X \mapsto \sum_{j=1}^{r_2} T_j X T_j^*, 
\end{align*}
where $K_j's$ and $T_j's$ are $k \times n$ matrices.

Equivalently, we can consider the following optimization problem: 
\begin{equation}\label{eq:qkdgeneric2}
    \begin{aligned}
    f(X) &= \Tr(\mathcal{L}_1(X) \ln(\mathcal{L}_1(X))) - \Tr(\mathcal{L}_1(X) \ln(\mathcal{L}_2(X))) \to \min,\\
    \Tr(A_i X) &= b_i,\,i=1,\ldots,m,\\
    X &\succeq 0. 
    \end{aligned}
\end{equation}

We experimentally apply our long-step path-following algorithm to \eqref{eq:qkdgeneric2} by solving the following auxiliary barrier family of optimization problems: for $\beta \geq 0$,
\begin{equation}\label{eq:qkdgenericbarrier}
    \begin{aligned}
    F_{\beta}(X) &= \beta f(X) - \ln\det(X)\to\min,\\
    \Tr(A_i X) &= b_i,\,i=1,\ldots,m,\\
    X &\succ 0. 
    \end{aligned}
\end{equation}
To guarantee the positive definiteness of $\mathcal{L}_1(X)$ and $\mathcal{L}_2(X)$, we use $\mathcal{L}_1(X) + \epsilon \cdot \id$ and $\mathcal{L}_2(X) + \epsilon \cdot \id$ instead, where $\epsilon$ is a very small positive number, e.g., $\epsilon = \num{1e-16}$. Note that in some computational settings, introduction of such perturbation may lead to certain instability. However, this is not very likely in our setting since first of all the perturbation is tiny, and secondly it is known that interior-point method is robust against small perturbation. Indeed, we do not see any instability in our numerical tests. 

We make the following conjecture based on our numerical results in \cref{table:qkd} and our previous results that the quantum relative entropy function is indeed compatible with the standard barrier $B(X) = -\ln\det(X)$ when either of the variables is fixed \cite{fayzhou19longstep, fayzhou19entanglement}.

\begin{conjecture}\label{conj:qkdselfcor}
$F_{\beta}(X)$ in \eqref{eq:qkdgenericbarrier} is a self-concordant function for each $\beta > 0$ and such self-concordance depends on the structure of the quantum relative entropy function and properties of the linear operators $\mathcal{L}_1$ and $\mathcal{L}_2$. 
\end{conjecture}

We strongly believe the conjecture is true and leave the theoretical proof for future work. The basic long-step path-following scheme remains the same as described in \cref{sec:longstepalg}, but the calculations of the gradient and Hessian of $F_{\beta}(X)$ are rather complicated, which we will present next. 

\subsection{Calculations of Gradient and Hessian}
Note that without loss of generality, we again only consider the real vector space of $n\times n$ symmetric matrices, i.e.,  $\mathbb{E} = \mathbb{S}^n$. The case of Hermitian matrices can be handled within the general Jordan algebraic scheme without much difficulty (see e.g. \cite{fayzhou19longstep}). Now let 
$$f_1(X) = \Tr(\mathcal{L}_1(X) \ln(\mathcal{L}_1(X))),\,f_2(X)= -\Tr(\mathcal{L}_1(X) \ln(\mathcal{L}_2(X))),$$

and 
$$B(X) = - \ln\det(X),$$
then
\begin{equation}
    \nabla F_{\beta}(X) = \beta (\nabla f_1(X) + \nabla f_2(X)) + \nabla B(X),
\end{equation}
and
\begin{equation}
        \Hessian_{F_{\beta}}(X) = \beta (\Hessian_{f_1}(X) + \Hessian_{f_2}(X)) + \Hessian_B(X).
\end{equation}
Next we show calculations for the three different components. Due to the importance of vectorization in implementation, we show its concrete forms along with the operator forms.

Let 
$$\mathcal{L}_1(X) = O_1\Lambda_1 O_1^T$$
be a spectral decomposition of $\mathcal{L}_1(X)$, where $\Lambda_1 = \diag(\lambda_1^{(1)},\ldots,\lambda_k^{(1)})$ and $O_1O_1^T = \id$. Similarly, let 
$$\mathcal{L}_2(X) = O_2\Lambda_2 O_2^T$$
be a spectral decomposition of $\mathcal{L}_2(X)$. Let $h(\lambda)=\ln(\lambda),\,\lambda>0$, and $h^{[1]}(\Lambda_i),\,i=1,2$, be the first divided difference introduced in \eqref{eq:1stdiv}. 

By \eqref{eq:vecid}, we have the vectorized forms of the linear operators $\mathcal{L}_1$ and $\mathcal{L}_2$:
\begin{equation}
    \begin{split}
        \widetilde{\mathcal{L}_1} &= \sum_{j=1}^{r_1} K_j \otimes \overline{K_j} = \sum_{j=1}^{r_1} K_j \otimes K_j,\\
        \widetilde{\mathcal{L}_2} &= \sum_{j=1}^{r_2} T_j \otimes \overline{T_j} = \sum_{j=1}^{r_2} T_j \otimes T_j,\\
    \end{split}
\end{equation}
where the complex conjugates are dropped in the second qualities because we only assume real matrices here. 

Consider the von Neumann entropy function
\begin{equation}\label{eq:trxlnx}
    f(X) = \Tr(X\ln(X)),\ X \succeq 0.
\end{equation}
Let $X = O\Lambda O^T$ be a spectral decomposition of $X$. Then the gradient of $f(X)$ is simply
\begin{equation}\label{eq:gradtrxlnx}
        \nabla f(X) = \id + \ln(X).
\end{equation}

For the Hessian, we have 
\begin{equation}
        \Hessian_f(X)(\xi,\xi) = \inp{\D\ln(X)(\xi)}{\xi},\,\forall \xi \in \mathbb{S}^n,
\end{equation}
and hence
\begin{equation}\label{eq:dlnx}
\begin{split}
    \Hessian_f(X)(\xi) &= \D\ln(X)(\xi)\\
                       &= O\left((O^T\xi O) \circ h^{[1]}(\Lambda) \right)O^T,\\
\end{split}
\end{equation}
where the last equality is shown in \eqref{eq:dlnxder} in \cref{sec:appendix}.

By \eqref{eq:gradgeneral} and \eqref{eq:gradtrxlnx}, we have the gradient of  $f_1(X)$:
\begin{equation}\label{eq:qkdgrad1}
        \nabla f_1(X) = \mathcal{L}_1^T\left(\id + \ln(\mathcal{L}_1(X))\right),
\end{equation}
and its vectorized form is simply
\begin{equation}
    \vecop(\nabla f_1(X)) = \widetilde{\mathcal{L}_1}^T \vecop\left(\id + \ln(\mathcal{L}_1(X))\right).
\end{equation}

By \eqref{eq:hessiangeneral1}, we have the Hessian of $\Hessian_{f_1}(X)$:
\begin{equation}\label{eq:qkdhession1}
    \begin{split}
            \Hessian_{f_1}(X)(\xi,\xi) &= \inp{\D\ln(\mathcal{L}_1(X)) (\mathcal{L}_1(\xi))}{\mathcal{L}_1(\xi)}\\
            & = \inp{\mathcal{L}_1^T \D\ln(\mathcal{L}_1(X)) (\mathcal{L}_1(\xi))}{\xi},\,\forall \xi \in \mathbb{S}^n,
    \end{split}
\end{equation}
and by \eqref{eq:dlnx}
\begin{equation}
\begin{split}
        \Hessian_{f_1}(X)(\xi) &= \mathcal{L}_1^T \D\ln(\mathcal{L}_1(X)) (\mathcal{L}_1(\xi))\\
        &= \mathcal{L}_1^T O_1\left((O_1^T \mathcal{L}_1(\xi) O_1) \circ h^{[1]}(\Lambda_1) \right)O_1^T.
\end{split}
\end{equation}

Using the vectorization procedure similar to \eqref{eq:gradmonotonevec} (take $C=\mathcal{L}_1(\xi)$), we have
\begin{equation}\label{eq:qkdf1hessianvecdev}
\begin{split}
        \vecop\left(\Hessian_{f_1}(X)(\xi)\right) 
        &= \widetilde{\mathcal{L}_1}^T \vecop\left(O_1\left((O_1^T \mathcal{L}_1(\xi) O_1) \circ h^{[1]}(\Lambda_1) \right)O_1^T\right)\\
        &= \widetilde{\mathcal{L}_1}^T(O_1\otimes O_1)\diag(\vecop(h^{[1]}(\Lambda_1)))(O_1\otimes O_1)^T\widetilde{\mathcal{L}_1}\vecop(\xi),
\end{split}
\end{equation}
and hence
\begin{equation}\label{eq:qkdf1hessianvec}
\Hessian_{f_1}(X) = \widetilde{\mathcal{L}_1}^T(O_1\otimes O_1)\diag(\vecop(h^{[1]}(\Lambda_1)))(O_1\otimes O_1)^T\widetilde{\mathcal{L}_1}.
\end{equation}

By \eqref{eq:gradgeneral} and product rule, we have 
\begin{equation}\label{eq:qkddf2}
    \D f_2(X)(\xi) = -\Tr\left[\mathcal{L}_1(\xi)\ln(\mathcal{L}_2(X))
    + \mathcal{L}_1(X) \D\ln(\mathcal{L}_2(X))(\mathcal{L}_2(\xi)) \right],\,\forall \xi \in \mathbb{S}^n,
\end{equation}
in which the first part is easy to take care of, and for the second part, use \eqref{eq:dtrclnx} (take $C=\mathcal{L}_1(X)$, replace $X$ with $\mathcal{L}_2(X)$ and $\xi$ with $\mathcal{L}_2(\xi)$) to get
\begin{equation}\label{eq:qkdgrad2}
    \nabla f_2(X) = -\mathcal{L}_1^T\ln(\mathcal{L}_2(X))-\mathcal{L}_2^T O_2\left((O_2^T\mathcal{L}_1(X)O_2)\circ h^{[1]}(\Lambda_2)\right)O_2^T.
\end{equation}
Similar to \eqref{eq:qkdf1hessianvecdev}, we have
\begin{equation}
\begin{split}
        \vecop(\nabla f_2(X))
        &= -\widetilde{\mathcal{L}_1}^T\!\!\vecop\left(\ln(\mathcal{L}_2(X))\right)-\widetilde{\mathcal{L}_2}^T\!\! \vecop\left(O_2\left((O_2^T\mathcal{L}_1(X)O_2)\circ h^{[1]}(\Lambda_2)\right) O_2^T\right)\\
        &= -\widetilde{\mathcal{L}_1}^T\!\!\vecop\left(\ln(\mathcal{L}_2(X))\right)-\widetilde{\mathcal{L}_2}^T\!\!(O_2 \otimes O_2) \diag(\vecop(h^{[1]}(\Lambda_2)))(O_2 \otimes O_2)^T \widetilde{\mathcal{L}_1}^T\!\!\vecop(X)
\end{split}
\end{equation}

By \eqref{eq:qkddf2} and product rule, we have
\begin{equation}
\begin{split}
    \Hessian_{f_2}(X)(\xi,\xi)
    & = \D^2 f_2(X)(\xi,\xi)\\
    & = -\Tr\left[\mathcal{L}_1(\xi) \D\ln(\mathcal{L}_2(X))(\mathcal{L}_2(\xi)) + \mathcal{L}_1(\xi) \D\ln(\mathcal{L}_2(X))(\mathcal{L}_2(\xi))\right.\\
    &\quad\quad\quad\quad+ \left.\mathcal{L}_1(X)\D^2\ln(\mathcal{L}_2(X))(\mathcal{L}_2(\xi), \mathcal{L}_2(\xi))\right]\\
    &= -\Tr\left[2\mathcal{L}_1(\xi) \D\ln(\mathcal{L}_2(X))(\mathcal{L}_2(\xi)) \right] - \Tr\left[\mathcal{L}_1(X)\D^2\ln(\mathcal{L}_2(X))(\mathcal{L}_2(\xi), \mathcal{L}_2(\xi))\right]\\
    &= -2\inp{\mathcal{L}_1^T \D\ln(\mathcal{L}_2(X))(\mathcal{L}_2(\xi))}{\xi} - \Tr\left[\mathcal{L}_1(X)\D^2\ln(\mathcal{L}_2(X))(\mathcal{L}_2(\xi), \mathcal{L}_2(\xi))\right].\\
\end{split}
\end{equation}
Let
\begin{align*}
    \Hessian_{f_{21}}(X)(\xi,\xi) &= -2\inp{\mathcal{L}_1^T \D\ln(\mathcal{L}_2(X))(\mathcal{L}_2(\xi))}{\xi},\text{ and }\\
    \Hessian_{f_{22}}(X)(\xi,\xi) &= - \Tr\left[\mathcal{L}_1(X)\D^2\ln(\mathcal{L}_2(X))(\mathcal{L}_2(\xi), \mathcal{L}_2(\xi))\right].
\end{align*}
By \eqref{eq:dlnx} we have 
\begin{equation}
\begin{split}
    \Hessian_{f_{21}}(X)(\xi) &= -2\mathcal{L}_1^T \D\ln(\mathcal{L}_2(X))(\mathcal{L}_2(\xi)),\\
    &= -2 \mathcal{L}_1^T O_2\left((O_2^T \mathcal{L}_2(\xi) O_2) \circ h^{[1]}(\Lambda_2)\right)O_2^T.
\end{split}
\end{equation}
Using symmetrization, and similar to the derivation in \eqref{eq:qkdf1hessianvecdev} and  \eqref{eq:qkdf1hessianvec}, we get
\begin{equation}\label{eq:qkdhessian21}
    \begin{split}
    \Hessian_{f_{21}}(X)
    &= -\,\widetilde{\mathcal{L}_1}^T(O_2\otimes O_2) \diag(\vecop(h^{[1]}(\Lambda_2))) (O_2\otimes O_2)^T \widetilde{\mathcal{L}_2}\\
    &\quad\, - \widetilde{\mathcal{L}_2}^T (O_2\otimes O_2) \diag(\vecop(h^{[1]}(\Lambda_2))) (O_2\otimes O_2)^T \widetilde{\mathcal{L}_1}.\\
    \end{split}
\end{equation}

By \eqref{eq:htrclnxder} (take $C = \mathcal{L}_1(X)$, replace $X$ with $\mathcal{L}_2(X)$ and $\xi$ with $\mathcal{L}_2(\xi)$), we get 
\begin{equation}
    \Hessian_{f_{22}}(X)(\xi,\xi) = \inp{O_2 \left(\intzerotoinfty{(D\tilde{C}D\tilde{\xi}D + D\tilde{\xi}D\tilde{C}D)}{t}\right) O_2^T}{\mathcal{L}_2(\xi)},\ \forall \xi \in \mathbb{S}^n,
\end{equation}
where $D = (\Lambda_2 + t\id)^{-1}$, $\tilde{C} = O_2^T\mathcal{L}_1(X)O_2$, and $\tilde{\xi} = O_2^T\mathcal{L}_2(\xi)O_2$.
Hence, 
\begin{equation}
    \Hessian_{f_{22}}(X)(\xi) = \mathcal{L}_2^T O_2 \left(\intzerotoinfty{(D\tilde{C}D\tilde{\xi}D + D\tilde{\xi}D\tilde{C}D)}{t}\right) O_2^T,\ \forall \xi \in \mathbb{S}^n.
\end{equation}
Using \eqref{eq:hessianmonotonevecdev}, we obtain
\begin{equation}
    \vecop\left(\Hessian_{f_{22}}(X)(\xi)\right) = \widetilde{\mathcal{L}_2}^T (O_2\otimes O_2) \intzerotoinfty{(D\tilde{C}D)\otimes D + D\otimes (D\tilde{C}D)}{t} (O_2\otimes O_2)^T \widetilde{\mathcal{L}_2}\vecop(\xi),
\end{equation}
hence
\begin{equation}\label{eq:qkdhessian22}
    \Hessian_{f_{22}}(X) = \widetilde{\mathcal{L}_2}^T(O_2\otimes O_2) \intzerotoinfty{(D\tilde{C}D)\otimes D + D\otimes (D\tilde{C}D)}{t}(O_2\otimes O_2)^T \widetilde{\mathcal{L}_2},
\end{equation}
where $D = (\Lambda_2 + t\id)^{-1}$ and $\tilde{C} = O_2^T\mathcal{L}_1(X)O_2$.

Combining \eqref{eq:qkdhessian21} and \eqref{eq:qkdhessian22}, we get
\begin{equation}
    \Hessian_{f_2}(X) = \Hessian_{f_{21}}(X) + \Hessian_{f_{22}}(X).
\end{equation}

As mentioned earlier (see \eqref{eq:hessianmonotonevec}), 
$$S = \intzerotoinfty{(D\tilde{C}D)\otimes D + D\otimes (D\tilde{C}D)}{t}$$
is a sparse block matrix with $(ij,kl)$-th entry: 
\begin{equation*}
    S_{ij,kl} = \delta_{kl}\tilde{C}_{ij}\Gamma_{ijl} + \delta_{ij}\tilde{C}_{kl}\Gamma_{jkl},
\end{equation*}
where 
$$\delta_{ij} = \left\{\begin{aligned}\,1\quad &\text{ if } i=j,\\ \,0\quad &\text{ if } i \neq j,
\end{aligned}
\right. \text{ and }\Gamma_{ijk} = -h^{[2]}(\lambda_i, \lambda_j, \lambda_k).$$
Lastly, we have 
\begin{equation}
    \begin{split}
        \vecop(\nabla B(X)) &= \vecop(-X^{-1}),\\
        \Hessian_B(X) &= X^{-1}\otimes X^{-1}.
    \end{split}
\end{equation}

\subsection{Numerical Results}
In this section, we present some of our numerical results. Again we use the analytic center as our initial point, which can be easily obtained by solving the following optimization problem (e.g., use SDPT3 \cite{sdpt3}): 

\begin{equation}\label{eq:miniacqkd}
\begin{aligned}
    f_{ac}(X) &= - \ln\det(X) \to \min,\\
    \inp{A_i}{X} &= b_i,\,i=1,\ldots,m,\\
    X &\succeq 0.
\end{aligned}
\end{equation}
Data are randomly generated \footnote{Data can be accessed here: \url{https://doi.org/10.13140/RG.2.2.19100.51847}} and without loss of generality, $\Tr(X) = 1$ is imposed (we assume that the feasible set is bounded). We set $\beta_0 = 0.1$, $\theta = 10$,  and $\epsilon = \num{1e-4}$ for all of our tests. Recall that the dimension of $X$ is $n \times n$. We use $k = 2n$ in our experiment. \Cref{table:qkd} shows numerical results for the QKD optimization problem \eqref{eq:qkdgeneric2} compared to the results obtained by using the \emph{quantum\_rel\_entr} function in \texttt{cvxquad} combined with SDP solver MOSEK \cite{mosek, fawzi2018cvxquad} (so far the most competitive approaches available for optimization problems involving quantum relative entropy). In \cref{table:qkd}, $nNewton$ is the number of total Newton steps used. $T_{ac}$ is the time for solving the analytic center and $T_{pf}$ is the time for running the long-step path-following algorithm. Both $T_{ac}$ and $T_{pf}$ are averaged over 20 repeated runs.

\begin{table}[tbhp]
\centering
\footnotesize{
\caption{Numerical results for QKD optimization problem \eqref{eq:minirelqkd}}\label{table:qkd}
\begin{tabular}{r r r r r r r c r r r }
\\\toprule
 & & & & &\multicolumn{3}{c}{Long-Step Path-Following} & \multicolumn{3}{r}{cvxquad $+$ mosek}\\

\midrule
$n$ & $k$ & $m$ & $r_1$ & $r_2$ & $T_{ac}$(s) & $T_{pf}$(s) & $nNewton$ &     $f_{min}$  &  Time(s) &     $f_{min}$ \\
\midrule

  4 & 8 & 2 & 2 & 2 & 0.15 & 0.03 & 6 & 0.2744 &  40.39 & 0.2744\\
  6 & 12 & 4 & 1 & 2 & 0.15 & 0.15 & 14 & 0.0498 & 2751.39 & 0.0498\\
  12 & 24 & 6 & 2 & 4 & 0.17 & 0.75 & 13 & 0.0440 & N/A & failed\\
  16 & 32 & 10 & 2 & 2 & 0.19 & 1.69 & 10 & 0.0511 & N/A & failed\\
  32 & 64 & 20 & 2 & 2 & 0.61 & 54.34 & 10 & 0.0332 & N/A & failed\\
\bottomrule
\end{tabular}
}
\end{table}

\begin{remarks}
\begin{enumerate}
    \item Our numerical results are quite stunning. For example, for the second test sample when $n=6$, our method is $9000$ times faster! For the larger dimensions, the other method simply can not solve the problem. Here \emph{failed} means that the program runs more than 10 hours without convergence. 
    \item During our numerical experiments, we notice that for most of the cases omitting the barrier term $B(X) = -\ln\det(X)$ does not affect the convergence. We suspect this is due to the fact that $-\Tr(C\ln(X)),\,X\succeq 0$,  is a self-concordant barrier for $C\succ 0$ (see \cite{fayzhou19entanglement}). However, the barrier term $B(X)$ does increase stability and accuracy of the algorithm.  
    \item Again our algorithm is implemented in Matlab, and all the numerical experiments are performed on a personal 15-in Macbook Pro with Intel core i7 and 16 GB memory.
\end{enumerate}
\end{remarks}

\section{Concluding Remarks}\label{sec:conclusion}
The difficulty of many optimization problems arising in quantum information theory stems from the fact that their objective functions are rather complicated nonlinear functions of several matrix arguments. In \cite{fayzhou19longstep, fayzhou19entanglement} and this paper we notice that many of these functions are compatible (in the sense of Nesterov and Nemirovskii) with the standard self-concordant barriers associated with symmetric cones. This observation, in principle, allows one to use structured interior-point algorithms for solving such optimization problems. To implement such algorithms one needs to be able to deal with very complicated Hessians. We show in detail how the analytic expressions of such Hessians can be derived along with their vectorized forms which are important for practical implementation. Although certain limitations on the size of the problem definitely exist, our extensive numerical experiments confirm that our long-step path-following algorithm is indeed robust and competitive. To the best of our knowledge our work is the first systematic attempt to use structured second-order methods for problems arising in quantum information theory. In comparison with first-order methods, our approach can solve comparable problems faster (asymptotic quadratic convergence) and with higher accuracy.

As for future work, first we would like to prove  \cref{conj:qkdselfcor}, which we strongly believe is true. Second, we formulate our algorithm in a very general setting (see \cref{thm:1,thm:2,thm:3}), but complexity estimates are available only for problems involving symmetric cones with standard self-concordant barrier functions. A natural question is whether it is possible to generalize the complexity estimates for the setting involving arbitrary self-concordant barriers. Last, it is of practical importance to explore size-reduction techniques, e.g., sparsity exploitation \cite{fujisawasparsity,vandenberghe2015chordal},  facial and symmetry reduction \cite{wolkowicz2017many,Bachoc2012}, to increase the problem size that can be realistically solved. For self-concordant functions like what we have, it would be particularly interesting to investigate the so-called \emph{Newton Sketch} \cite{pilanci2015newton}, which can substantially reduce the computational cost by performing approximate Newton steps.

\section*{Acknowledgement}
This research is supported in part by Simmons Foundation Grant 275013.

\bibliographystyle{abbrvnat}
\bibliography{refs}

\appendix 
\section{}\label{sec:appendix}
For $C \succeq 0$, let
\begin{equation}
    f(X) = - \Tr(C\ln(X)),\ X\succeq 0.
\end{equation}
We will show how to derive the gradient and Hessian of $f(X)$ by using the integral representation of $\ln(X)$.  

Let $\id$ be the identity matrix of the same size as $X$. Since 
\begin{equation*}
    \ln(x) = \intzerotoinfty{\frac{1}{1+t} - (x+t)^{-1}}{t},
\end{equation*}
we have, for $X \succeq 0$, the integral representation of $\ln(X)$:
\begin{equation}
            \ln(X) = \intzerotoinfty{\left[\frac{1}{1+t}\id - (X+t\id)^{-1}\right]}{t},
\end{equation}
with which we can derive its first and second Fr\'{e}chet derivatives: 
\begin{align}
        \D \ln (X) (\xi) &= \intzerotoinfty{(X+t\id)^{-1}\xi(X+t\id)^{-1}}{t},\label{eq:1stdevln}\\
        \D^2 \ln (X)(\xi,\xi) &= -2\intzerotoinfty{(X+t\id)^{-1}\xi(X+t\id)^{-1}\xi(X+t\id)^{-1}}{t}.\label{eq:2nddevln}
\end{align}

Consider a spectral decomposition of $X$, 
$$X = U\Lambda U^T,$$
where $\Lambda = \diag(\lambda_1,\ldots,\lambda_n)$ and $UU^T = \id$. Then 
\begin{equation*}
    (X + t\id)^{-1} = UDU^T,
\end{equation*}
where $D = (\Lambda + t\id)^{-1}$. For $C \succeq 0$, let 
$$\tilde{C} = U^T C U,\ \tilde{\xi} = U^T \xi U,$$ and $$d_i = (\lambda_i + t)^{-1},\,i=1,\ldots,n.$$

Furthermore, let $h(\lambda) = \ln(\lambda),\,\lambda > 0$, and denote the first divided difference matrix by $h^{[1]}(\Lambda)$, where
\begin{equation*}
\begin{split}
    [h^{[1]}(\Lambda)]_{ij} &= h^{[1]}(\lambda_i,\lambda_j)\\
    &= \left\{\begin{aligned}
             \frac{h(\lambda_j) - h(\lambda_i)}{\lambda_j - \lambda_i}&, &  &\lambda_j \ne \lambda_i, \\
             h'(\lambda_i)&, &  &\lambda_j = \lambda_i.
        \end{aligned}\right.
\end{split}
\end{equation*}

We have 
\begin{equation}\label{eq:intdxid}
    \intzerotoinfty{D\tilde{\xi} D}{t} = h^{[1]}(\Lambda) \circ \tilde{\xi},
\end{equation}
where $\circ$ is the Schur product: for $m \times n$ matrices $A$ and $B$,
$$[A \circ B]_{ij} = A_{ij}B_{ij}.$$
Indeed, 
\begin{equation*}
    \begin{split}
        \intzerotoinfty{[D\tilde{\xi} D]_{ij}}{t} &= \intzerotoinfty{d_i\tilde{\xi}_{ij}d_j}{t}\\
        &= \tilde{\xi}_{ij}\intzerotoinfty{d_i d_j}{t}\\
        &= \tilde{\xi}_{ij}\intzerotoinfty{(\lambda_i + t)^{-1}(\lambda_j + t)^{-1}}{t}\\
        &= \tilde{\xi}_{ij} \cdot \left\{\begin{aligned}
             \frac{\ln(\lambda_j) - \ln(\lambda_i)}{\lambda_j - \lambda_i}&, &  &\lambda_j \ne \lambda_i \\
             \frac{1}{\lambda_j}&, &  &\lambda_j = \lambda_i.
        \end{aligned}\right.\\
        &= \tilde{\xi}_{ij}\cdot h^{[1]}(\lambda_i,\lambda_j).
    \end{split}
\end{equation*}

Then
\begin{equation}\label{eq:dtrclnx}
    \begin{split}
        \D f(X)(\xi) 
        &= - \Tr(C\D\ln(X)(\xi))\\
        &\!\!\!\overset{\eqref{eq:1stdevln}}{=}\!\!\! - \Tr\left(C\intzerotoinfty{(X+tI)^{-1}\xi(X+tI)^{-1}}{t}\right)\\
        &= - \Tr \left(\intzerotoinfty{CUDU^T \xi UDU^T}{t}\right)\\
        &= - \intzerotoinfty{\Tr(U^TCUDU^T \xi UD)}{t}\\
        &= - \intzerotoinfty{\Tr(\tilde{C}D\tilde{\xi}D)}{t}\\
        &= - \Tr\left(\tilde{C}\intzerotoinfty{D\tilde{\xi}D}{t}\right)\\
        &\!\!\!\overset{\eqref{eq:intdxid}}{=}\!\!\! -\Tr\left(\tilde{C}\left(h^{[1]}(\Lambda) \circ \tilde{\xi}\right)\right)\\
        &= - \Tr\left(\left(\tilde{C} \circ h^{[1]}(\Lambda)\right) \tilde{\xi}\right)\\
        &= \inp{-U\left(\tilde{C} \circ h^{[1]}(\Lambda)\right)U^T}{\xi},
    \end{split}
\end{equation}
where in the second last equality we used properties of Schur product (\cite[p.~306]{horn_johnson_1991}):
$$\Tr((A\circ B)C^T) = \Tr((A\circ C)B^T).$$
In a similar spirit, we can also easily get
\begin{equation}\label{eq:dlnxder}
\begin{split}
        \D \ln(X)(\xi) &= \intzerotoinfty{UDU^T \xi UDU^T}{t}\\
        &= U \left(\intzerotoinfty{DU^T \xi UD}{t}\right) U^T\\
        &\!\!\!\overset{\eqref{eq:intdxid}}{=}\!\! U\left(h^{[1]}(\Lambda) \circ (U^T \xi U)\right) U^T
\end{split}
\end{equation}

Now for the Hessian of $f(X)$, we have
\begin{equation}\label{eq:htrclnxder}
    \begin{split}
        \D^2 f(X)(\xi,\xi) 
        &= - \Tr\left(C\D^2 \ln(X)(\xi,\xi)\right)\\
        &\!\!\!\overset{\eqref{eq:2nddevln}}{=} \!\!\!-\Tr\left(C\left(-2\intzerotoinfty{(X+t\id)^{-1}\xi(X+t\id)^{-1}\xi(X+t\id)^{-1}}{t}\right)\right)\\
        &= - \Tr\left(C\left(-2 \intzerotoinfty{UDU^T\xi UDU^T\xi UDU^T}{t}\right)\right)\\
        &= 2 \intzerotoinfty{\Tr(\tilde{C}D\tilde{\xi}D\tilde{\xi}D)}{t}\\
        &= 2 \intzerotoinfty{\Tr(D\tilde{C}D\tilde{\xi}D\tilde{\xi})}{t}\\
        &= 2 \intzerotoinfty{\Tr\left(\frac{D\tilde{C}D\tilde{\xi}D + D\tilde{\xi}D\tilde{C}D}{2}\cdot\tilde{\xi}\right)}{t}\\
        &= \Tr\left(\left(\intzerotoinfty{(D\tilde{C}D\tilde{\xi}D + D\tilde{\xi}D\tilde{C}D)}{t}\right)\tilde{\xi}\right)\\
        &= \inp{U \left(\intzerotoinfty{(D\tilde{C}D\tilde{\xi}D + D\tilde{\xi}D\tilde{C}D)}{t}\right) U^T}{\xi}\\
        &= \inp{\Hessian_f(X)(\xi)}{\xi}.
    \end{split}
\end{equation}
Hence,
\begin{equation}\label{eq:htrclnx}
    \Hessian_f(X)(\xi) = U \left(\intzerotoinfty{(D\tilde{C}D\tilde{\xi}D + D\tilde{\xi}D\tilde{C}D)}{t}\right) U^T.
\end{equation}

\end{document}